\documentclass[reqno,11pt]{amsart} 
\usepackage{amsmath}
\usepackage{amssymb}
\usepackage{amsthm}
\usepackage{verbatim}
\usepackage{xcolor}
\usepackage{graphics}

\newcommand\NoBlackBoxes{\global\overfullrule0pt}

\NoBlackBoxes
\theoremstyle{plain}

\begin{document}

\title{ENERGY BOUNDS FOR KANTOROVICH TRANSPORT \\
DISTANCES WITH CONVEX COST FUNCTIONS
}

\author{Sergey G. Bobkov$^{1,3}$}
\thanks{1) 
School of Mathematics, University of Minnesota, Minneapolis, MN, USA,
bobkov@math.umn.edu. 
}

\author{Friedrich G\"otze$^{2,3}$}
\thanks{2) Faculty of Mathematics,
Bielefeld University, Germany,
goetze@math-uni.bielefeld.de.
}

\thanks{3) Research supported by the GRF – SFB 1283/2 2021 – 317210226}


\begin{abstract}
Energy bounds for Kantorovich transport distances are developed for convex cost 
functions. The main results extend estimates due to M. Ledoux for the Kantorovich 
distances $W_p$.
\end{abstract}

\maketitle
\markboth{Sergey G. Bobkov and Friedrich G\"otze}{Energy bounds}

\def\theequation{\thesection.\arabic{equation}}
\def\E{{\mathbb E}}
\def\R{{\mathbb R}}
\def\C{{\mathbb C}}
\def\P{{\mathfrak P}}
\def\Z{{\mathbb Z}}
\def\L{{\mathbb L}}
\def\T{{\mathcal T}}

\def\G{\Gamma}

\def\Ent{{\rm Ent}}
\def\var{{\rm Var}}
\def\Var{{\rm Var}}

\def\H{{\rm H}}
\def\Im{{\rm Im}}
\def\Tr{{\rm Tr}}
\def\s{{\mathfrak s}}

\def\k{{\kappa}}
\def\M{{\cal M}}
\def\Var{{\rm Var}}
\def\Ent{{\rm Ent}}
\def\O{{\rm Osc}_\mu}

\def\ep{\varepsilon}
\def\phi{\varphi}
\def\vp{\varphi}
\def\F{{\cal F}}

\def\be{\begin{equation}}
\def\en{\end{equation}}
\def\bee{\begin{eqnarray*}}
\def\ene{\end{eqnarray*}}

\thispagestyle{empty}

\section{{\bf Introduction}}
\setcounter{equation}{0}

\vskip2mm
\noindent
Let $\P(\R^n)$ denote the collection of all (Borel) probability measures on the Euclidean 
space $\R^n$ and let $\Pi(\mu,\nu)$ denote the collection of all probability measures $\pi$ 
on $\R^n \times \R^n$ with given marginals $\mu$ and $\nu$ in $\P(\R^n)$, 
called couplings or transference (transport) plans. The Kantorovich transport problem
amounts to describe or estimate the transportation cost
$$
\T_c(\mu,\nu) = \inf_{\pi \in \Pi(\mu,\nu)}\int\!\!\!\int  c(x,y)\,d\pi(x,y),
$$
where $c:\R^n \times \R^n \rightarrow [0,\infty)$ is a fixed Borel measurable function, 
called cost function.

Here cost functions of power type, that is $c_p(x,y) = |x-y|^p$, $p \geq 1$, give rise to the
transport metrics $W_p = \T_{c_p}^{1/p}$ in the space of all probability measures 
$\mu$ on $\R^n$ with finite $p$-th absolute moment $\int |x|^p\,d\mu(x)$.
As it turns out, $W_p$ may be estimated in terms of the dual Sobolev norm
\be
\|\nu - \mu\|_{H^{-1,p}(\mu)} = \sup\Big\{\int f\,d\nu - \int f\,d\mu:
\|\nabla f\|_{L^q(\mu)} \leq 1\Big\},
\en
where the supremum runs over all 
$C^\infty$-smooth, compactly supported functions $f$ on $\R^n$ such that
$$
\|\nabla f\|_{L^q(\mu)} = \Big(\int |\nabla f|^q\, d\mu\Big)^{1/q} \leq 1, \quad
q = p/(p-1)
$$
(with usual understanding in the case $p=1$). Recently M. Ledoux has proposed the following 
remarkable bound on the power Kantorovich distance $W_p$ called energy estimate 
(motivated by the bounds for convergence rates of empirical measures with respect to
$W_2$).

\vskip5mm
{\bf Theorem 1.1} (\cite{L}). {\sl For all $\mu$ and $\nu$ in $\P(\R^n)$ with
finite $p$-th absolute moment, $p \geq 1$,
\be
W_p(\mu,\nu) \leq p\,\|\nu - \mu\|_{H^{-1,p}(\mu)}.
\en
}

\vskip0mm
This bound was actually derived in the setting of (complete connected) Riemannian manifolds $M$
without boundary, assuming that 
$\nu$ is absolutely continuous with respect to $\mu$ on $M$ (a condition that may be removed).
Putting $g = \varphi - 1$, $\varphi = d\nu/d\mu$, M. Ledoux defined the  dual Sobolev norm
equivalently on functions $g$ on $M$ with $\mu$-mean zero by
$$
\| g \|_{H^{-1,p}(\mu)} = \Big(\int_M |\nabla (-T^{-1} g)|^p\,d\mu\Big)^{1/p}.
$$
Here $T$ represents a second order linear differential operator which serves as the infinitesimal 
generator of the Markov semi-group associated to $\mu$, defined by
$$
\int_M u (-T v)\,d\mu = \int_M \left<\nabla u,\nabla v\right> d\mu,
$$
and $T^{-1}$ denotes the inverse of the non-negative operator $-T$.

The aim of this note is to develop the ``infimum-convolution" approach used in \cite{L} and
as a result -- to extend (1.2) to more general cost functions of the form $c(x,y) = L(x-y)$, 
that is, for the optimal transportation cost
$$
\T_L(\mu,\nu) \, = \inf_{\pi \in \Pi(\mu,\nu)} \int\!\!\!\int  L(x-y)\,d\pi(x,y).
$$
We assume that $L$ is a convex function on $\R^n$ such that $L(0)=0$, $L(x)>0$ for 
$x \neq 0$, and $\frac{L(x)}{|x} \rightarrow \infty$ as $|x| \rightarrow \infty$. 
Note that the transport distance $\T_L(\mu,\nu)$ does not need to be symmetric in $(\mu,\nu)$.

It is also assumed that $L$ satisfies a two-sided $\Delta_2$-condition $L(2x) \leq cL(x)$,
valid for all $x \in \R^n$ with some constant $c$. In this and only this case the associated 
Young function
$$
\Phi_L(r) \, = \sup_{x \neq 0} \, \frac{L(rx)}{L(x)},\quad r \geq 0,
$$
is finite. Necessarily it is convex on the positive half-axis (as the supremum of convex functions).
In this setting the definition (1.1) has a natural generalization
\be
\|\nu - \mu\|_{H^{-1,L}(\mu)} = \sup\Big\{\int f\,d\nu - \int f\,d\mu:
\int L^*(\nabla f)\,d\mu \leq 1\Big\},
\en
where
$$
L^*(x) = \sup_{y \in \R^n} \big[\left<x,y\right> - L(y)\big], \quad x \in \R^n,
$$
is the Legendre transform of $L$. The integrals in (1.3) are well-defined and finite 
(since $f$'s are bounded), although in general the supremum may be infinite. Thus, 
in terms of the Luxemburg pseudo-norm generated by $L^*$, the value 
$c = \|\nu - \mu\|_{H^{-1,L}(\mu)}$ represents an optimal constant in the inequality
$$
\int f\,d\nu - \int f\,d\mu \leq c\,\|\nabla f\|_{L^*(\mu)}
$$
in the class of all $C^\infty$-smooth functions $f$ on $\R^n$ with compact support.

First we consider the ``absolutely continuous" case.

\vskip5mm
{\bf Theorem 1.2.} {\sl Let $\mu$ and $\nu$ in $\P(\R^n)$ satisfy
$\int\!\!\int L(x-y)\,d\mu(x)\,d\nu(y) < \infty$. If $\nu$ is absolutely continuous with respect to 
$\mu$ which in turn is absolutely continuous with respect to the Lebesgue measure on $\R^n$, 
then
\be
\T_L(\mu,\nu) \leq A\, \Phi_L\big(\|\nu - \mu\|_{H^{-1,L}(\mu)}\big)
\en
with constant $A = \Phi_L(\Phi_L'(1+))$.
}

\vskip5mm
In particular, if $L(x) = \|x\|^p$ for a norm $\|\cdot\|$ on $\R^n$, we have $\Phi_L(r) = r^p$
and $A = p^p$. In this case, (1.4) includes the Ledoux' bound (1.2). 

Using a smoothing argument, one can give a slightly weakened version of Theorem 1.2, 
by removing the condition of absolute continuity, even without weakening of the constant
for the homogeneous case of order $p > 1$ where $L(rx) = r^p L(x)$ for all $x \in \R^n$ and $r \geq 0$.

Introduce the one-sided derivatives $p_+ = \Phi'_L(1+)$, $p_- = \Phi'_L(1-)$
and the associated quantity
$$
\gamma = \sup \big\{a + b: a+br \leq \min(r^{p_+},r^{p_-}) \ {\rm for \ all} \ r \geq 0\big\}.
$$
In particular, $0 < \gamma \leq 1$, and $\gamma = 1$ if and only if $p_+ = p_- = p$.
The latter may be shown to be equivalent to the property that $L$ is positive homogeneous 
of order $p$.
Note that since $\Phi_L(1)=1$, we have $p_+ \geq p_- \geq 1$. If $p_- > 1$, then 
necessarily $L$ satisfies
the super-linear growth condition $L(x)/|x| \rightarrow \infty$ as $|x| \rightarrow \infty$.

\vskip5mm
{\bf Theorem 1.3.} {\sl If $p_- > 1$, then for all probability measures $\mu$ and $\nu$ 
on $\R^n$ such that $\int\!\!\int L(x-y)\,d\mu(x)\,d\nu(y) < \infty$, the inequality $(1.4)$
holds true with constant
$$
A = \Phi_L(p_+)\, \Phi_L(1/\gamma).
$$
}

\vskip0mm
The paper is a organized as follows. First we recall basic results on Kantorovich duality,
infimum-convolution operators and Hamilton-Jacobi equations (Sections 2-6).
The two-sided $\Delta_2$-condition is discussed  in Sections 7-8, including one special
boundary value problem about the Young functions $\Phi_L$. Basic definitions and 
results about the Luxemburg and Orlicz pseudo-norms for vectors valued-functions
are collected in  Section 9. Sections 10-11 are devoted to properties of dual Sobolev
norms (the lower semi-continuity and continuity along convolutions
needed in the smoothing argument). Proofs of Theorems 1.2 and 1.3 are completed in
Section 12.

\vskip7mm
\section{{\bf Kantorovich Duality}}
\setcounter{equation}{0}

\vskip2mm
\noindent
For reader's convenience, here we collect several general facts about the Kantorovich duality.
Given a separable metric space $M$, denote by $\P(M)$ the collection of all Borel probability 
measures $\mu$ on $M$. The product space $M \times M'$ of two separable metric spaces
is metrizable and separable. Fix a lower semi-continuous (cost) function 
$c:M \times M' \rightarrow [0,\infty)$.

Given $\mu \in \P(M)$ and $\nu \in \P(M')$, denote by $\Pi(\mu,\nu)$ the collection 
of all Borel probability measures $\pi$ on $M \times M'$ with marginals $\mu$ and $\nu$. 
Like in the Euclidean space, the problem is how
to describe or estimate the transportation cost
\be
\T_c(\mu,\nu) = \inf_{\pi \in \Pi(\mu,\nu)}\int\!\!\!\int  c(x,y)\,d\pi(x,y).
\en

\vskip2mm
The set $\Pi(\mu,\nu)$ is non-empty, since it contains, for example, the product measure 
$\pi = \mu \otimes \nu$, which represents a transference transport plan. Hence
$\T_c(\mu,\nu)$ is well-defined. If both $M$ and $M'$ are complete separable metric spaces,
$\Pi(\mu,\nu)$ is compact in $\P(M \times M')$ in the weak topology, by
Prokhorov's compactness criterion. Since the cost function is lower semi-continuous, 
the infimum in (2.1) is attained for some plan $\pi$, called an optimal plan.

The next theorem is known as the Kantorovich duality, cf. Villani \cite{V1}, Theorem 1.3.

\vskip5mm
{\bf Proposition 2.1.} {\sl Let $M,M'$ be complete separable metric spaces.
For all $\mu \in \P(M)$ and $\nu \in \P(M')$,
\be
\T_c(\mu,\nu) = \sup \bigg[\int_M f\,d\mu + \int_{M'} g\,d\nu\bigg],
\en
where the supremum is running over all functions $f \in L^1(\mu)$ and
$g \in L^1(\nu)$ such that
\be
f(x) + g(y) \leq c(x,y)
\en
for $\mu$-almost all $ x \in M$ and $\nu$-almost all $y \in M'$. Subject to $(2.3)$,
the supremum in $(2.2)$ may be restricted to functions $f,g$ that are bounded 
and continuous.
}

\vskip5mm
Applying (2.3), we obtain that, for any transference transport plan
 $\pi$,
\bee
\int_M f\,d\mu + \int_{M'} g\,d\nu
 & = &
\int_M \int_{M'} (f(x) + g(y))\,d\pi(x,y) \\
 & \leq &
\int_M \int_{M'} c(x,y)\,d\pi(x,y).
\ene
Taking the supremum over all admissible pairs $(f,g)$ and then the infimum
over all $\pi$, we arrive at
$$
\sup \bigg[\int_M f\,d\mu + \int_{M'} g\,d\nu\bigg] \leq \T_c(\mu,\nu).
$$

In order to reverse this inequality, note that the optimal
function $g$ for a fixed function $f$  and the optimal
function $f$ for a fixed function $g$ in (2.3) are given by
\bee
g(y)
 & = &
Q f(y) \, = \, \inf_{x \in M} \big[c(x,y) - f(x)\big], \\
f(x)
 & = &
Q g(x) \, = \, \inf_{y \in M'} \big[c(x,y) - g(x)\big]. \\
\ene
The next assertion is taken from \cite{A-G-S}, cf. Theorem 6.1.5 and Remark  6.1.6.

\vskip5mm
{\bf Proposition 2.2.} {\sl Let $M,M'$ be separable metric spaces, and let the cost function
satisfy
$$
\int\!\!\!\int  c(x,y)\,d\mu(x)\,\nu(y) < \infty
$$
for given measures $\mu \in \P(M)$ and $\nu \in \P(M')$. There exists
a maximizing pair $(f,Qf)$ in $(2.2)$, and if $\pi$ is an optimal plan, then
\be
f(x) + Qf(y) = c(x,y) \quad \pi-{\sl a.e.} \ in \ M \times M'. 
\en
Moreover, if there exists a function $f \in L^1(\mu)$ such that $(2.4)$
holds for a transference transport  plan $\pi$, then this plan is optimal.
}

\vskip5mm
The question of uniqueness of optimal plans was also addressed in the literature.
The following assertion holds for the setting of Euclidean spaces
$M = M' = \R^n$, cf. \cite{A-G-S}, Th.\,6.2.4.

\vskip5mm
{\bf Proposition 2.3.} {\sl Assume that the cost function has the form
$c(x,y) = L(x-y)$ for some strictly convex function $L:\R^n \rightarrow [0,\infty)$.
Suppose furthermore that
$$
\int\!\!\!\int  L(x-y)\,d\mu(x)\,\nu(y) < \infty, \quad \mu, \nu \in \P(\R^n).
$$
If $\mu$ is absolutely 
continuous with respect to the Lebesgue measure on $\R^n$, there exists 
a unique plan $\pi \in \Pi(\mu,\nu)$ minimizing the double integral in $(2.1)$.
Moreover, there is a Borel map $T:\R^n \rightarrow \R^n$ such that 
$\pi$ is the push forward of $\mu$ under the map $x \rightarrow (x,T(x))$.
}

\vskip5mm
In particular, $\nu$ is the push forward of $\mu$ under the map $T$.
In the sequel, we will be interested in the distances
$\T_L(\mu,\nu)$ defined for the cost functions $c(x,y) = L(x-y)$.

\vskip7mm
\section{{\bf Infimum-Convolution Operators}}
\setcounter{equation}{0}

\vskip2mm
\noindent
Let $L:\R^n \rightarrow [0,\infty)$ be a convex function such that $L(0) = 0$.
We associate with $L$ the infimum-convolution operators 
\be
Q_t f(x) = \inf_{y \in \R^n} \bigg[f(y) + t L\Big(\frac{x-y}{t}\Big)\bigg], \quad x \in \R^n,
\en
where $t>0$ is a parameter. They are well-defined for all functions $f:\R^n \rightarrow \R$. 
At the initial time $t=0$, one puts $Q_0 f = f$ suggested by continuity properties of $Q_t$ 
which we discuss later.

In general, $Q_t f \leq f$. Since $Q_t f$ represents the infimum of continuous functions
$x \rightarrow f(y) + t L(\frac{x-y}{t})$, it is lower semi-continuous, hence, Borel measurable. 
Moreover, if 
$f$ is bounded, then $Q_t f$ is bounded as well and has finite Lipschitz semi-norm, so it is continuous
(this will be shown below).

According to Proposition 2.1 (Kantorovich duality) and the particular case (2.2),
$$
\T_L(\mu,\nu) = \sup \bigg[\int g\,d\nu - \int f\,d\mu\bigg],
$$
where the supremum is running over all functions $f \in L^1(\mu)$ and $g \in L^1(\nu)$,
or equivalently, over all bounded continuous $f$ and $g$
such that $g(x) - f(y) \leq L(x-y)$ for all $x,y \in \R^n$. 
Given $f$, here the optimal choice is $g = Q_1 f(x)$. Thus, we obtain an equivalent
description of the Kantorovich duality in terms of the infimum-convolution operators:

\vskip5mm
{\bf Proposition 3.1.} {\sl For all $\mu$ and $\nu$ in $\P(\R^n)$, we have
\be
\T_L(\mu,\nu) = \sup \bigg[\int Q_1 f\,d\nu - \int f\,d\mu\bigg],
\en
where the supremum is running over all bounded continuous functions $f$ on $\R^n$.
}

\vskip5mm
By the convexity of $L$,
$$
\inf_{y \in \R^n} \Big[t L\Big(\frac{x-y}{t}\Big) + s L\Big(\frac{y-z}{s}\Big)\Big] =
(t+s)\, L\Big(\frac{x-z}{t+s}\Big), \quad x,z \in \R^n, \ t,s>0,
$$
which ensures that these operators form a semigroup: $Q_t Q_s f = Q_{t + s} f$.
Under the assumption $\frac{L(x)}{|x|} \rightarrow \infty$ as $|x| \rightarrow \infty$,
this semigroup has a generator, which may be shown to be
\be
\frac{d}{dt}\,Q_t f \Big|_{t=0} = -L^*(\nabla f)
\en
for all $f$ from a suitable domain.
Note that, by the growth condition on $L$, the Legendre transform $L^*$ is convex and 
finite on the whole space $\R^n$. In addition, it is non-negative and satisfies
$\frac{L^*(x)}{|x|} \rightarrow \infty$ as $|x| \rightarrow \infty$ (cf. e.g. \cite{E}).

The equality (3.3) represents the heart of the Hamilton-Jacobi equation. 
In order to give more precise statements, we need to look at basic properties of the
infimum-convolution operators. For example, the function $t \rightarrow Q_t f(x)$
is non-increasing, and $|Q_t f(x)| \leq M$ as long as $|f(x)| \leq M$ for all $x \in \R^n$ 
with some $M \geq 0$. In this case the infimum in (3.1) may be restricted to certain balls 
in $\R^n$, as long as $L(x) \rightarrow \infty$ as $|x| \rightarrow \infty$. 
Define
\be
R_L(r) = \sup\big\{|x|: L(x) \leq r\big\}, \quad r \geq 0,
\en
which is the radius of the smallest closed ball in $\R^n$ with center at the origin
containing the convex compact set $K_L(r) = \{x \in \R^n: L(x) \leq r\}$.

\vskip4mm
{\bf Proposition 3.2.} {\sl Suppose that $L(x) \rightarrow \infty$ as $|x| \rightarrow \infty$. 
If $|f(x)| \leq M$ for all $x \in \R^n$ with some constant $M \geq 0$, then for any $t>0$, 
\be
Q_t f(x) = \inf_{|x-y| \leq r} \Big[f(y) + t L\Big(\frac{x-y}{t}\Big)\Big], \quad 
r = t R_L(2M/t).
\en
}

\vskip0mm
{\bf Proof.} It should be clear that  $R_L(r)$ is non-decreasing and tends to infinity as 
$r \rightarrow \infty$. Putting $z = x-y$, we need to choose $r \geq 0$ such that 
$|z| > r \Rightarrow t L(z/t) > 2M$. Equivalently, $L(z/t) \leq 2M/t \Rightarrow |z| \leq r$, or
$z \in t K_L(2M/t) \Rightarrow |z| \leq r$. By the definition (3.4), the optimal value is
$r = t R_L(2M/t)$, in which case the expression in the square brackets in (3.5) will be larger than 
$-M + 2M = M$. Since $Q_t f(x) \leq M$, the points $y$ such that $|x-y| > r$ may be therefore
excluded from the infimum in (3.1).
\qed

\vskip5mm
As a consequence of this observation, we shall show:

\vskip5mm
{\bf Proposition 3.3.} {\sl Suppose that $L(x) \rightarrow \infty$ as $|x| \rightarrow \infty$. 
For any $t>0$, the supremum
\be
\sup \bigg[\int Q_t f\,d\nu - \int f\,d\mu\bigg]
\en
within the class of all bounded continuous functions coincides with the supremum taken over all 
$C^\infty$-smooth non-negative functions $f$ on $\R^n$ with compact support.
}

\vskip5mm
Since the difference of integrals in (3.6) does not change when adding any 
constant to $f$, one may require that $f \geq 0$. 

First, let us see that the supremum in (3.6) may be restriced to the class of all bounded continuous functions 
$f \geq 0$ on $\R^n$ with compact support. If $f \geq 0$ is bounded and continuous, 
consider the function $f_R = f g_R$, where $g_R:\R^n \rightarrow [0,1]$ is continuous, supported
on the ball $B_{R+1}= \{x \in \R^n: |x| \leq R+1\} $, and such that $g_R = 1$ on $B_R$. 
Then $f_R$ will be continuous and supported on $B_{R+1}$. In addition, $f_R = f$ on $B_R$ 
and, due to (3.5), $Q_t f_R = Q_t f$ on $B_{R-r}$ for $R>r = tR_L(2M/t)$. As a consequence,
$$
\int Q_t f_R\,d\nu - \int f_R\, d\mu \rightarrow \int Q_t f\,d\nu - \int f\, d\mu\quad 
(R \rightarrow \infty).
$$

Next, let $f \geq 0$ on $\R^n$ be bounded, continuous with compact support. 
Take a $C^\infty$-smooth function
$\omega:\R^n \rightarrow [0,\infty)$ supported on the unit ball $B_1$
and such that $\int \omega(x)\,dx = 1$. Then the functions $\omega_\ep(x) = \ep^n \omega(\ep x)$
with parameter $\ep \in (0,1]$ are also $C^\infty$-smooth and are supported on the balls $B_\ep$. 
Consider the convolutions
\be
f_\ep(x) = (f * \omega_\ep)(x) = \int f(x-y)\,\omega_\ep(y)\,dy = \int \omega_\ep(x-y)\,f(y)\,dy,
\en
which are bounded by $M = \sup_x f(x)$, $C^\infty$-smooth, and  have a compact support. 
In addition, $f_\ep(x) \rightarrow f(x)$ as $\ep \rightarrow 0$ for all $x \in \R$.
Hence, by the dominated convergence theorem,
$$
\lim_{\ep \rightarrow 0}\int f_\ep\, d\mu = \int f d\mu.
$$

We need a similar property for $Q_t f_\ep$. Since $Q_t f \leq M$ and $Q_t f_\ep \leq M$,
$$
\Big|\int Q_t f(x)\,d\nu(x) - \int_{|x| \leq R} Q_t f(x)\,d\nu\Big| \leq 
M \nu\{|x| > R\} \rightarrow 0 \quad (R \rightarrow \infty),
$$
and the same inequality is true for $Q_t f_\ep$. Hence
\be
\Big|\int (Q_t f_\ep - Q_t f)\,d\nu\Big| \leq \sup_{x \in B_R} |Q_t f_\ep(x) - Q_t f(x)| +
2M \nu\{|x| > R\}.
\en

Using the modulus of continuity
$$
\delta_R(\ep) = \sup\{|f(x) - f(y)|: x\in B_R, \ |x-y| \leq \ep\},
$$
it follows from (3.7) that $|f_\ep(x) - f(x)| \leq \delta_R(\ep)$ for all $x \in B_R$.
Applying (3.5) to $f$ and $f_\ep$, we therefore obtain that
$$
|Q_t f_\ep(x)  - Q_t f(x)| \leq \sup_{|h| \leq r} |f_\ep(x+h) - f(x+h)| \leq \delta_{R+r}(\ep)
$$
for all $x \in B_R$. Applying this in (3.8), we arrive at
$$
\Big|\int (Q_t f_\ep - Q_t f)\,d\nu\Big| \leq \delta_{R+r}(\ep) + 2M \nu\{|x| > R\}.
$$
Here the right-hand side can be made arbitrary small using the uniform 
continuity of $f$ on every ball $B_R$. Hence, we may conclude that
$\int Q_t f_\ep\,d\nu \rightarrow \int Q_t f\,d\nu$ as $\ep \rightarrow 0$, and
$$
\int Q_t f_\ep\,d\nu - \int f_\ep\, d\mu \rightarrow \int Q_t f\,d\nu - \int f\, d\mu.
$$
Thus, the supremum in (3.6) may be restricted to the functions of the form $f_\ep$.
\qed

\vskip7mm
\section{{\bf Basic Properties of Infimum-Convolution Operators}}
\setcounter{equation}{0}

\vskip2mm
\noindent
Here we collect some general properties of infimum-convolution operators.
Let us start with monotonicity with respect to the time variable. Below in this section, 
we assume that $L$ is a non-negative convex function on $\R^n$ such that $L(0)=0$.

\vskip5mm
{\bf Proposition 4.1.} {\sl Given $f:\R^n \rightarrow \R$ and a point $x \in \R^n$, the function
$u(t) = Q_t f(x)$ is non-increasing in $t>0$. 
}

\vskip2mm
{\bf Proof.}
By the convexity, $L$ admits the representation
$L(h) = \sup_{(a,b) \in \mathfrak F} (\left<a,h\right> + b)$
for some set $\mathfrak F \subset \R^n \times \R$, so,
$$
t L(h/t) = \sup_{(a,b) \in \mathfrak F} (\left<a,h\right> + bt).
$$
Since $L(0)=0$, necessarily $b\leq 0$ for any couple $(a,b) \in \mathfrak F$.
This shows that this  function is non-increasing in $t$ as the supremum of non-increasing
functions. Hence, $u(t)$  is non-increasing being the infimum of non-increasing
functions $t \rightarrow f(x-h) + t L(h/t)$.
\qed

\vskip5mm
A function $f$ on $\R^n$ will be called Lipschitz, if it has a finite Lipschitz semi-norm 
$\|f\|_{\rm Lip}$ with respect to the Euclidean distance. In this case, it follows
 from the definition (3.1) that
\bee
Q_t f(x) 
 & = &
\inf_{h \in \R^n} \Big[f(x-h) + t L\Big(\frac{h}{t}\Big)\Big] \\
 & \geq &
\inf_{h \in \R^n} \Big[f(x) - \|f\|_{\rm Lip}\,|h| + t L\Big(\frac{h}{t}\Big)\Big] \\
 & = &
f(x) - t \sup_{y \in \R^n} \Big[\|f\|_{\rm Lip}\,|y| - L(y)\Big].
\ene
Using $|y| = \sup_{|\theta|=1} \left<\theta,y\right>$,
the last expression may be related to the Legendre transform of~$L$.

\vskip5mm
{\bf Proposition 4.2.} {\sl Given a Lipschitz function $f$ on $\R^n$, for all $x \in \R^n$ and $t>0$,
\be
Q_t f(x) \geq f(x) - t \, \sup_{|\theta|=1} L^*\big(\|f\|_{\rm Lip}\,\theta\big).
\en
}

We will need the following fact, which is recorded in Evans \cite{E}, Lemma 2 on page 127.

\vskip5mm
{\bf Proposition 4.3.} {\sl If $f$ is a Lipschitz function on $\R^n$, then so is $Q_t f$, and moreover, 
\be
\|Q_t f\|_{\rm Lip} \leq \|f\|_{\rm Lip} \quad {\sl for \ all} \ \, t>0.
\en
In addition,
the function $u_x(t) = Q_t f(x)$ is Lipschitz on the half-axis $[0,\infty)$ with 
\be
\|u_x\|_{\rm Lip} \leq \sup_{|\theta|=1} L^*(\|f\|_{\rm Lip}\,\theta) \quad {\sl for \ all} \ \, x \in \R^n,
\en
provided that the last supremum is finite.
}

\vskip5mm
{\bf Corollary 4.4.} {\sl If $f$ is Lipschitz on $\R^n$ and the supremum in $(4.3)$ is 
finite, then the function $u(x,t) = Q_t f(x)$ is differentiable almost everywhere on 
$\R^n \times (0,\infty)$. Also, for any $x \in \R^n$, the function 
$u_x(t) = Q_t f(x)$ is locally absolutely continuous. Moreover, the limit
$$
u_x'(t) = \limsup_{\ep \downarrow 0} \frac{u_x(t+\ep) - u_x(t)}{\ep}
$$
exists for almost all $t>0$, is bounded by $\sup_{|\theta|=1} L^*(\|f\|_{\rm Lip}\,\theta)$ 
in absolute value, and may serve as a Radon-Nikodym derivative of $u_x(t)$.
}

\vskip5mm
{\bf Proof of Proposition 4.3.}
Since $f$ may decay at infinity not faster than a linear function, the infimum in the definition
(3.1) is attained at some point $y = y_t(x)$. Then, for any $x' \in \R^n$, 
\bee
Q_t f(x') - Q_t f(x)
 & = &
\inf_z \Big[f(z) + t L\Big(\frac{x'-z}{t}\Big)\Big] - \Big[f(y) + t L\Big(\frac{x-y}{t}\Big)\Big] \\
 & \leq &
f(x'-x+y) - f(y) \, \leq \, \|f\|_{\rm Lip}\,|x'-x|,
\ene
where we chose $z = x'-x+y$ in the last infimum (to equalize $x'-z$ and $x-y$).
Interchanging the roles of $x$ and $x'$, we arrive at (4.2).

For the second claim, we apply (4.1) together with (4.2) to get that, for any $\ep>0$,
\bee
Q_{t+\ep} f(x) \, = \, Q_\ep Q_t f(x) 
 & \geq &
Q_t f(x) - \ep\, \sup_{|\theta|=1} L^*\big(\|Q_t f\|_{\rm Lip}\,\theta\big) \\
 & \geq &
Q_t f(x) - \ep\, \sup_{|\theta|=1} L^*\big(\|f\|_{\rm Lip}\,\theta\big).
\ene
The resulting inequality is equivalent to (4.3).
\qed

\vskip7mm
\section{{\bf Behavior of $Q_t$ at the Initial Time $t=0$ and its Generator}}
\setcounter{equation}{0}

\vskip2mm
\noindent
In this section, we consider the behaviour of $Q_t f(x)$ for small values of $t>0$. 
Assume that $L$ is a non-negative convex function on $\R^n$ such that $L(0)=0$
and $L(x)/|x| \rightarrow \infty$ as $|x| \rightarrow \infty$. 

Define
$$
f_-(x) = \liminf_{y \rightarrow x} f(y) = \lim_{\ep \rightarrow 0} \inf_{|x-y| < \ep} f(y).
$$

\vskip3mm
{\bf Proposition 5.1.} {\sl Given a bounded function $f$ on $\R^n$, for any point $x \in \R^n$,
\be
\lim_{\ep \downarrow 0}\, Q_\ep f(x) = f_-(x).
\en
Thus, $Q_\ep f(x) \rightarrow f(x)$ as $\ep \downarrow 0$ if and only if $f$ is lower
semi-continuous at $x$.
}

\vskip5mm
{\bf Proof.} First note that, for every $\ep>0$,
$$
Q_\ep f(x) \leq \liminf_{y \rightarrow x} \bigg[f(y) + \ep L\Big(\frac{x-y}{\ep}\Big)\bigg] = f_-(x),
$$
so, $\limsup_{\ep \downarrow 0} Q_\ep f(x) \leq f_-(x)$. 
For the opposite direction, define  $f_\delta(x) = \inf_{|x-y| \leq \delta} f(y)$. If $|x-y| \leq \delta$,
$$
f(y) + \ep L\Big(\frac{x-y}{\ep}\Big) \geq f(y) \geq f_\delta(x).
$$
If $|x-y| > \delta$, then
$$
f(y) + \ep L\Big(\frac{x-y}{\ep}\Big) \geq m + \ep\,T_L(\delta/\ep),
$$
where $m = \inf_x f(x)$ and $T_L(r) = \inf_{|x| > r} L(x)$, $r \geq 0$.
Both inequalities yield
\be
Q_\ep f(x) \geq \min\big\{f_\delta(x), m + \ep\,T_L(\delta/\ep)\big\}.
\en
By the growth condition on $L$, we have $T_L(r)/r \rightarrow \infty$ as $r \rightarrow \infty$, 
so that $\ep\,T_L(\delta/\ep) \rightarrow \infty$ as $\ep \rightarrow 0$ for any fixed $\delta>0$. 
Hence, from (5.2),
$$
\liminf_{\ep \downarrow 0}\, Q_\ep f(x) \geq f_\delta(x).
$$ 
It remains to let $\delta \rightarrow 0$.
\qed

\vskip5mm
{\bf Proposition 5.2.} {\sl 
If a bounded function $f$ on $\R^n$ is differentiable at the point $x \in \R^n$, then
\be
\lim_{\ep \downarrow 0} \frac{Q_\ep f(x) - f(x)}{\ep} = -L^*(\nabla f(x)).
\en
}

\vskip2mm
{\bf Proof.} Let $M = \sup_x |f(x)|$. Fix $x \in \R^n$. The assumption on $f$ ensures that
\be
f(x-h) = f(x) - \left<\nabla f(x),h\right> + \delta(h)\,|h|
\en
with $\delta(h) \rightarrow 0$ as $h \rightarrow 0$. By the definition of $Q_\ep$, 
for all $h \in \R^n$,
$$
Q_\ep f(x) \leq f(x-h) + \ep L(h/\ep).
$$
Hence, putting $h_\ep = \ep u$, $u \in \R^n$, we have
$$
f(x) - Q_\ep f(x) 
 \geq
f(x) - f(x-h_\ep) - \ep L(u),
$$
which gives, by (5.4),
\bee
\liminf_{\ep \downarrow 0}\, \frac{f(x) - Q_\ep f(x)}{\ep} 
 & \geq & 
\lim_{\ep \downarrow 0} \frac{f(x) - f(x - h_\ep)}{\ep} - L(u) \\
 & = &
\left<\nabla f(x),u\right> - L(u).
\ene
The left-hand side does not depend on $u$. So,
taking the supremum over all $u$ on the right-hand side, we arrive at
$$
\liminf_{\ep \downarrow 0}\, \frac{f(x) - Q_\ep f(x)}{\ep} \geq L^*(\nabla f(x)).
$$

For the opposite direction, we use the growth assumption on $L$ which yields
$$
\ep(r) = \sup_{|x| \geq r} \frac{|x|}{L(x)} \downarrow 0 \quad {\rm as} \ r \rightarrow \infty.
$$
As $L(x) \geq |x|/\ep(|x|)$, we have $L(x) \leq r \Rightarrow |x|/\ep(|x|) \leq r$. 
The latter implies that $|x| = o(r)$ uniformly over all $x$ such that $L(x) \leq r$, which yields 
\be
\frac{1}{r}\,R_L(r) = \frac{1}{r}\, \sup\big\{|x|: L(x) \leq r\big\} \rightarrow 0 \quad {\rm as} \ 
r \rightarrow \infty.
\en

We now apply Proposition 3.2. Recall that, by (3.5), for any $\ep>0$,
\be
Q_\ep f(x) \, = \inf_{|h| \leq r(\ep)} \Big[f(x-h) + \ep L\Big(\frac{h}{\ep}\Big)\Big], \quad 
r(\ep) = \ep R_L(2M/\ep).
\en
By (5.5), $r(\ep) \rightarrow 0$ as $\ep \rightarrow 0$.
For $0 < |h| \leq r(\ep)$, put $\theta = h/|h|$ and $v = \nabla f(x) - \delta(h) \theta$.
Applying (5.4) and the definition of $L^*$, we get
\bee
f(x) - f(x-h) 
 & = &
\left<\nabla f(x) - \delta(h) \theta,h\right> \\
 & = & 
\ep \left<v,h/\ep\right> \, \leq \, \ep L^*(v) + \ep L(h/\ep).
\ene
By (5.6), this gives
\bee
\frac{f(x) - Q_\ep f(x)}{\ep} 
 & = &
\sup_{|h| \leq r(\ep)} \Big[\frac{f(x) - f(x-h)}{\ep} - \ep L(h/\ep)\Big] \\
 & \leq &
\sup_{|h| \leq r(\ep)}  L^*\big(\nabla f(x) - \delta(h) \theta\big) \, \rightarrow \, 
L^*(\nabla f(x))
\ene
as $\ep \rightarrow 0$, where we used $\sup_{|h| \leq r(\ep)} |\delta(h)| \rightarrow 0$
together with continuity of $L^*$. Hence
$$
\limsup_{\ep \downarrow 0}\, \frac{f(x) - Q_\ep f(x)}{\ep} \leq L^*(\nabla f(x)).
$$
\qed

\vskip7mm
\section{{\bf Hamilton-Jacobi Equations}}
\setcounter{equation}{0}

\vskip2mm
\noindent
As in Proposition 5.2, suppose that $L$ is a non-negative convex function on $\R^n$ 
such that $L(0)=0$ and $L(x)/|x| \rightarrow \infty$ as $|x| \rightarrow \infty$.

If $f$ is a Lipschitz function on $\R^n$, then, according to
Proposition 4.2, for any $x \in \R^n$, the function 
$$
u(t) = Q_t f(x) = \inf_{h \in \R^n} \bigg[f(x-h) + t L\Big(\frac{h}{t}\Big)\bigg]
$$
is Lipschitz and therefore absolutely continuous on the positive half-axis $[0,\infty)$. In addition, 
it is non-increasing. As a consequence, for any $x \in \R^n$, the partial derivative
$$
\frac{\partial}{\partial t}\, Q_t f(x) = \lim_{\ep \rightarrow 0}
\frac{Q_{t+\ep} f(x) - Q_t f(x)}{\ep}
$$
exists and is finite for almost all $t > 0$.

\vskip5mm
{\bf Proposition 6.1.} {\sl Let $f$ be a bounded Lipschitz function on $\R^n$. 
For all points $(x,t) \in \R^n \times (0,\infty)$ of differentiability of $u(x,t) = Q_t f(x)$ 
$($hence for almost all $(x,t)$ \!\!$)$, 
\be
\frac{\partial}{\partial t}\, Q_t f(x) = -L^*\big(\nabla Q_t f(x)\big).
\en
}

\vskip0mm
{\bf Proof.} By Proposition 5.2, applied to the function $Q_t f$ in place of $f$,
$$
\lim_{\ep \downarrow 0} \frac{Q_{t+\ep} f(x) - Q_t f(x)}{\ep} \, = \,
\lim_{\ep \downarrow 0} \frac{Q_\ep Q_t f(x) - Q_t f(x)}{\ep} \, = \,  -L^*\big(\nabla Q_t f(x)\big).
$$
Here the first limit exists by the assumption that $u(x,t)$ is differentiable with respect
to the time variable at a fixed point $t>0$.
\qed

\vskip5mm
The equality in (6.1) represents the Hamilton–Jacobi equations 
with initial condition
\[
\left\{\begin{array}{ll}
u_t = - L(\nabla u) & {\rm on} \ \ \R^n \times (0,\infty), \\
u = f & {\rm on} \ \ \R^n \times \{0\},
\end{array}\right.\]
where $f(x)$ is a given function on $\R^n$ and $u = u(x,t)$ is an unknown function on
$\R^n \times [0,\infty)$, with the usual notations $u_t = \partial_t u$,
$\nabla u = \nabla_x u$ (the gradient with respect to the space variable). The solution to this 
initial value problem) may therefore be given in the form of the infimum-convolution
$u(x,t) = Q_t f(x)$ which is also called the Hopf–Lax formula.
Here the boundedness assumption on $f$ may be removed (cf. Evans \cite{E}, pp. 128-129).

On the other hand, for some classes of convex functions $L$, the Lipschitz condition, as well as the
differentiability assumption may also be removed (still keeping boundedness of $f$), by a proper
definition of the gradient function $\nabla Q_t f(x)$. Some refinements of Proposition 6.1
for the quadratic function $L(x) = |x|^2$ can be found in \cite{B-C-G}.

One important application of the Hamilton–Jacobi equation (6.1) is the following
assertion about intepolation of functions and integrals.

\vskip5mm
{\bf Corollary 6.2.} {\sl Let $f$ be a bounded Lipschitz function on $\R^n$ and let
$\nu$ be an absolutely continuous probability measure on $\R^n$. For any $t>0$,
\be
\int (Q_t f - f)\,d\nu = -\int_0^t \int_{\R^n} L^*\big(\nabla Q_s f(x)\big)\,d\nu(x)\,ds.
\en
Moreover, for almost all $x \in \R^n$,
\be
Q_t f(x) - f(x) = -\int_0^t L^*\big(\nabla Q_s f(x)\big)\,ds.
\en
}

{\bf Proof.} The integrand in (6.2)-(6.3) is well-defined for almost all $(x,s)$ in view of the
a.e. differentiabilty of the function $(x,s) \rightarrow Q_s f(x)$ and due to the absolute continuity
of the measure $\nu$. At all such points $(x,s)$ we also have $|\nabla Q_s f(x)| \leq \|f\|_{\rm Lip}$,
by Proposition 4.3. So, the integrand in (6.2)-(6.3) represents a Borel measurable bounded function
defined on a set in $\R^n \times (0,\infty)$ of full Lebesgue measure.

Now, consider the function 
$$
\psi(t) = \int (Q_t f(x) - f(x))\,d\nu(x), \quad t \geq 0.
$$ 
We have $\psi(0) = 0$ and, by Proposition 4.3, 
$$
|\psi(t) - \psi(s)| \leq c\,|t-s|, \quad c = \sup_{|\theta|=1} L^*(\|f\|_{\rm Lip}\,\theta), \ \ t,s \geq 0.
$$
That is, $\psi$ is Lipschitz and hence locally absolutely continuous, with representation
\be
\psi(t) = \int_0^t \psi'(s)\,ds,
\en
where $\psi'$ is a Radon-Nikodym derivative of $\psi$, which coincides with the usual derivative
of $\psi$ existing and finite for almost all $t>0$, and satisfying $|\psi'(t)| \leq c$.

In order to identify the usual derivative, fix $t>0$, $\ep>0$, and write
$$
\frac{\psi(t+\ep)-\psi(t)}{\ep} = \int \frac{Q_{t+\ep} f(x)  -Q_t f(x)}{\ep}\,d\nu(x).
$$
By Proposition 6.1, for a set $\Omega \subset \R^n \times (0,\infty)$ of full Lebesgue
measure, as $\ep \downarrow 0$,
$$
g_{\ep,t}(x) = \frac{Q_{t+\ep} f(x) - Q_t f(x)}{\ep} \rightarrow -L^*\big(\nabla Q_t f(x)\big), \quad
(x,t) \in \Omega.
$$
In addition, by Proposition 4.3, $|g_{\ep,t}(x)|$ is bounded by the constant $c$.
The section $\Omega(t) = \{x \in \R^n: (x,t) \in \Omega\}$ has a full Lebesgue measure
for almost all $t>0$. Hence, we may apply the Lebesgue dominated convergence 
theorem and conclude that, for almost every $t>0$, the function $\psi$ has a finite derivative
at this point given by
$$
\psi'(t) = -\int_{\R^n} L^*\big(\nabla Q_t f(x)\big)\,d\nu(x).
$$
It remains to apply the identity (6.4), and then we arrive at (6.2).

In particular, for any $x \in \R^n$ and $\ep>0$, (6.2) yields
$$
\frac{1}{|B(x,\ep)|} \int_{B(x,\ep)} (Q_t f(x) - f(x))\,dx = -\frac{1}{|B(x,\ep)|} 
\int_{B(x,\ep)} \bigg[\int_0^t L^*\big(\nabla Q_s f(x)\big)\,ds\bigg]\,dx,
$$
where $B(x,\ep)$ is the Euclidean ball with center at $x$ and radius $\ep$, and where
$|B(x,\ep)|$ denotes its $n$-dimensional volume. By the Lebesgue differentiation theorem,
for almost all $x$, both normalized integrals are convergent as $\ep \downarrow 0$ to the
corresponding integrands, and then we obtain (6.3).
\qed

\vskip7mm
\section{{\bf The Two-Sided $\Delta_2$-Condition and Associated Young Functions}}
\setcounter{equation}{0}

\vskip2mm
\noindent
Let $L$ be a convex function on $\R^n$ such that $L(0)=0$ and $L(x)>0$ for $x \neq 0$. 
As was already defined, we say that it satisfies a two-sided $\Delta_2$-condition, if
\be
L(2x) \leq cL(x), \quad x \in \R^n,
\en
with some constant $c$. In the usual definition of the $\Delta_2$-condition, (7.1) is required
to hold for all sufficiently large $x$, i.e., when $|x| \geq x_0$ for some $x_0>0$.
In that case, $L$ may be redefined near zero in such way that a new function would satisfy (7.1)
for all $x$, and then it will generate the same Orlicz space (cf. \cite{K-R} for the one dimensional
setting). Due to the current requirement, the behavior of the function $L$ can be controlled
not only at infinity, but also near zero (viewed as a second side).

In this section we recall some results related to the two-sided $\Delta_2$-condition which
have been recently discussed in \cite{B-G}. This condition is equivalent to the
property that the function
$$
\Phi_L(r) \, = \sup_{x \neq 0} \, \frac{L(rx)}{L(x)},\quad r \geq 0,
$$
is finite on the positive half-axis. We call it the Young function associated to $L$.
Since the two-sided $\Delta_2$-condition is needed in derivation of energy estimates 
for the transport distances, let us list basic properties of $\Phi_L$.

\vskip5mm
{\bf Proposition 7.1.} {\sl  Under the two-sided $\Delta_2$-condition the following properties 
hold true for the associated Young function $\Phi = \Phi_L$:

\vskip2mm
$a)$ \ $\Phi$ is non-negative, non-decreasing, convex, and satisfies $\Phi(0)=0$, $\Phi(r)>0$ for $r > 0$;

$b)$ \ $\Phi(1)=1$ and $1 \leq \Phi'(1-) \leq \Phi'(1+)$;

$c)$ \ $\Phi$ is sub-multiplicative: $\Phi(rs) \leq \Phi(r) \Phi(s)$ for all $r,s \geq 0$;

$d)$ \ $\Phi$ has a sub-polynomial growth: putting $p_- = \Phi'(1-)$ and $p_+ = \Phi'(1+)$, we have
\be
\Phi(r) \leq r^{p_-} \ (0 \leq r\leq 1),  \qquad \Phi(r) \leq r^{p_+} \ (r \geq 1).
\en
}

As a closely related function, introduce
$$
\Psi_L(r) \, = \inf_{x \neq 0} \, \frac{L(rx)}{L(x)} = \frac{1}{\Phi_L(1/r)}, \quad r>0.
$$
It possesses similar properties.

\vskip3mm
{\bf Proposition 7.2.} {\sl  Under the two-sided $\Delta_2$-condition the following properties 
hold true for the function $\Psi = \Psi_L$:

\vskip2mm
$a)$ \ $\Psi$ is non-negative, non-decreasing, and satisfies $\Psi(0)=0$, $\Psi(r)>0$ for $r > 0$;

$b)$ \ $\Psi(1)=1$ and $\Psi'(1-) = \Phi'(1+)$, $\Psi'(1+) = \Phi'(1-)$;

$c)$ \ $\Psi$ is super-multiplicative: $\Psi(rs) \geq \Psi(r) \Psi(s)$ for all $r,s \geq 0$;

$d)$ \ $\Psi$ has a super-polynomial growth: 
\be
\Psi(r) \geq r^{p_+} \ (0 \leq r\leq 1),  \qquad \Psi(r) \geq r^{p_-} \ (r \geq 1).
\en
}

\vskip0mm
In contrast with $\Phi_L$, the function $\Psi_L$ does not need be convex.

Since $\Psi(r) \leq \Phi(r)$, the comparison of (7.2) with (7.3) leads to the following conclusion.

\vskip5mm
{\bf Corollary 7.3.} {\sl  Under the two-sided $\Delta_2$-condition, the function $\Phi_L$
is differentiable at the point $r=1$ and has derivative $p = \Phi_L'(1)$, if and only if 
$\Phi_L(r) = \Psi_L(r) = r^p$ for all $r \geq 0$, that is, if $L$ is positive homogeneous of order $p$:
\be
L(rx) = r^p L(x), \quad x \in \R^n, \ r \geq 0.
\en
}

Necessarily $p = p_- = p_+$.
Typical examples of such functions are given by $L(x) = \|x\|^p$, $x \in \R^n$,
with parameter $p \geq 1$, where $\|\cdot\|$ is an arbitrary norm on $\R^n$. 
Non-homogeneous examples include the functions of the form $L(x) = V(\|x\|)$,
where $V$ is a Young function, i.e. $V:[0,\infty) \rightarrow [0,\infty)$ is convex, $V(0) = 0$, 
$V(r) > 0$ for $r>0$. Then $\Phi_L(r) = \sup_{s > 0} \, \frac{V(rs)}{V(s)}$.

The two-sided $\Delta_2$-condition may be characterized in terms of the directional
derivatives
$$
L'(x,y) = \lim_{\ep \downarrow 0} \frac{L(x+\ep y) - L(x)}{\ep}, \quad x,y \in \R^n,
$$
which are well-defined and finite due to the convexity of $L$.
Namely, (7.1) is equivalent to
\be
\sup_{x \neq 0} \, \frac{L'(x,x)}{L(x)} \, < \, \infty.
\en
In this case, this supremum does not exceed $\Phi_L'(1+)$. If $L$ is differentiable,
(7.5) simplifies to
\be
\sup_{x \neq 0} \, \frac{\left<\nabla L(x),x\right>}{L(x)} < \infty.
\en

Restricting these suprmuma to regions $|x| > x_0 > 0$, (7.5)--(7.6) will provide the
characterization of the usual $\Delta_2$-condition. In dimension $n=1$, and when $L$ is even, 
we obtain a well-known characterization of the $\Delta_2$-condition:
$\sup_{x > x_0} L'(x) x/L(x) < \infty$,
where $L'(x)$ denotes the right derivative of $L$ at the point $x$.

Another question we need to address is when the Legendre transform $L^*$
satisfies the two-sided $\Delta_2$-condition. Recall that 
$1 \leq p_- \leq p_+$ where $p_- = \Phi_L'(1-)$ and $p_+ = \Phi_L'(1+)$.

\vskip5mm
{\bf Proposition 7.4.} {\sl If $L$ satisfies the two-sided $\Delta_2$-condition with $p_- > 1$,
$L^*$ satisfies this condition as well. In this case
\be
\Phi_{L^*}(1+) \leq \frac{p_-}{p_- - 1}, \quad \Phi_{L^*}(1-) \leq \frac{p_+}{p_+ - 1}.
\en
}

This statement can be refined for the class of homogeneous convex functions.

\vskip5mm
{\bf Proposition 7.5.} {\sl If $L$ is positive homogeneous of order $p \geq 1$, it
satisfies the two-sided $\Delta_2$-condition. In this case, $L^*$ satisfies the two-sided $\Delta_2$-condition
if and only if $p > 1$, and then it is positive homogeneous of order $q = p/(p-1)$.
}

\vskip7mm
\section{{\bf A Boundary Value Problem under the Two-Sided $\Delta_2$-Condition}}
\setcounter{equation}{0}

\vskip2mm
\noindent
Let us return to the Young function $\Phi = \Phi_L$ associated to the convex function 
$L$ on $\R^n$ such that $L(0)=0$ and $L(x)>0$ for $x \neq 0$. Assuming the two-sided
$\Delta_2$-condition, one immediate consequence of Proposition 7.1 is the following:

\vskip5mm
{\bf Corollary 8.1.} {\sl The function $\Phi$ has an inverse
$\Phi^{-1}:[0,\infty) \rightarrow [0,\infty)$, which is a concave increasing function with
$\Phi^{-1}(0) = 0$, $\Phi^{-1}(1) = 1$ and $\Phi^{-1}(\infty) = \infty$. In addition,
\be
\int_\delta^\infty \frac{ds}{s\, \Phi^{-1}(s)} < \infty \ \ {\sl for \ all} \ \delta>0, \qquad
\int_0^\infty \frac{ds}{s\, \Phi^{-1}(s)} = \infty.
\en
Hence, for any $c>0$ there is a unique solution $\delta>0$ to the equation
\be
\int_\delta^\infty \frac{ds}{s\, \Phi^{-1}(s)} = \frac{1}{c}.
\en
}

\vskip0mm
Indeed, by $d)$, we have $\Phi^{-1}(s) \geq s^{1/p}$ for all $s \geq 1$ with $p = \Phi'(1+)$, 
so that the first integral in (8.1) is convergent. Since $\Phi^{-1}(s) \leq 1$
on the interval $0 \leq s \leq 1$, the second integral is divergent.

We use this observation to establish a specific boundary value problem, which is needed
to establish energy-type bounds.

\vskip5mm
{\bf Proposition 8.2.} {\sl Given $c>0$, let $\delta>0$ be the solution to the equation $(8.2)$.
There exists a unique increasing smooth function $\theta:[0,1] \rightarrow [0,1]$ such that
$\theta(0) = 0$, $\theta(1) = 1$, and
\be
(1 - \theta(t))\,\Phi\Big(\frac{c\theta'(t)}{1 - \theta(t)}\Big) = \delta, \quad 0 \leq t < 1.
\en
The derivative $\theta'(t)$ is bounded on $[0,1)$.
}

\vskip4mm
For example, for the power functions $\Phi(r) = r^p$, $p \geq 1$, the solution to (8.2) 
is given by $\delta = (cp)^p$, and the equation (8.3) is solved explicitly by 
$\theta(t) = 1 - (1-t)^p$.

\vskip5mm
{\bf Proof.} (8.3) may be rewritten as the equation
\be
\theta'(t) = U(\theta(t)), \quad 0 \leq t < 1,
\en
where
$$
U(y) = \frac{1}{c}\,(1-y)\,\Phi^{-1}\Big(\frac{\delta}{1-y}\Big), \quad 0 \leq y < 1.
$$
Since, by the concavity of $\Phi^{-1}$, the function $\Phi^{-1}(s)/s$ is non-increasing 
in $s > 0$, $s \Phi^{-1}(1/s)$ is non-decreasing, so, the function $U(y)$ is positive and 
non-increasing on $[0,1)$. In particular, $U(y) \leq U(0) = \frac{1}{c}\,\Phi^{-1}(\delta)$,
and by (8.4), $\theta'(t) \leq \frac{1}{c}\,\Phi^{-1}(\delta)$, proving the last claim.

Consider the positive increasing smooth function
$$
R(z) = \int_0^z \frac{1}{U(y)}\,dy, \quad 0 \leq z \leq 1,
$$
and let $R^{-1}$ be its inverse. We have $R(\theta(t))'  = \theta'(t)/U(\theta(t)$, so that
(8.4)  may be stated as the equation $R(\theta(t))' = 1$ on
$[0,1)$. But, since $R(0) = \theta(0) = 0$,  this is equivalent to
$$
R(\theta(t)) = t \ \Longleftrightarrow \ \theta(t) = R^{-1}(t).
$$
Thus, $\theta = R^{-1}$ provides the solution to the equation (8.4).
By construction, $\theta(0) = R^{-1}(0) = R(0) = 0$. The other requirement
$\theta(1)=1$ is equivalent to $R(1) = 1$. And indeed,
$$
R(1) \, = \, \int_0^1 \frac{1}{U(y)}\,dy \, = \, 
\int_0^1 \frac{c}{(1-y)\,\Phi^{-1}(\frac{\delta}{1-y})}\,dy \, = \,
c\int_\delta^\infty \frac{1}{s\,\Phi^{-1}(s)}\,ds.
$$
The latter expression is equal to 1, if and only if (8.2) is satisfied (which was assumed).
\qed

\vskip5mm
We now turn to the question of how to estimate $\delta = \delta(c)$ in (8.2). As a function 
of $c$, it is increasing and continuous with $\delta(0) = 0$ and $\delta(\infty) = \infty$.

\vskip5mm
{\bf Proposition 8.3.} {\sl Given $c>0$, the solution $\delta = \delta(c)$ in $(8.2)$ satisfies
\be
\delta \leq A \Phi(c), \quad A = \Phi(\Phi'(1+)).
\en
}

\vskip0mm
Note that there is an equality in (8.5) for the power functions $\Phi(r) = r^p$, $p \geq 1$,
in which case $A = p^p$.

\vskip2mm
{\bf Proof.} By Proposition 7.1, the inverse function $\Phi^{-1}$ is super-multiplicative, so that
$\Phi^{-1}(\delta t) \geq \Phi^{-1}(\delta) \, \Phi^{-1}(t)$ for all $t \geq 0$. Hence, changing 
the variable $s = \delta t$
and using $\Phi^{-1}(t) \geq t^{1/p}$ for $t \geq 1$ with $p = \Phi'(1+)$, we have
\bee
\int_\delta^\infty \frac{ds}{s\, \Phi^{-1}(s)} 
 & = &
\int_1^\infty \frac{dt}{t\, \Phi^{-1}(\delta t)} \\
 &  \leq &
\frac{1}{\Phi^{-1}(\delta)} \int_1^\infty \frac{dt}{t\, \Phi^{-1}(t)}
 \ \leq \
\frac{1}{\Phi^{-1}(\delta)} \int_1^\infty \frac{dt}{t^{1 + 1/p}} \ = \ \frac{p}{\Phi^{-1}(\delta)}.
\ene
As a consequence, from (8.2) it follows that $\frac{1}{c} \leq \frac{p}{\Phi^{-1}(\delta)}$.
In terms of $\delta$, this inequality is solved as $\delta \leq \Phi(pc) \leq \Phi(p)\, \Phi(c)$, 
where we used the sub-multiplicativity of $\Phi$.
\qed

\vskip7mm
\section{{\bf Luxemburg and Orlicz Pseudo-Norms for Vector-Valued Functions}}
\setcounter{equation}{0}

\vskip2mm
\noindent
Let $L:\R^n \rightarrow [0,\infty]$ be a lower semi-continuous convex function such that 
$L(0)=0$, $L(x)>0$ for all $x \neq 0$, which is finite in some neighborhood of the origin.
The conjugate transform
$$
L^*(x) = \sup_{y \in \R^n} \big[\left<x,y\right> - L(y)\big], \quad x \in \R^n,
$$
is also non-negative, convex, but may vanish near zero and take an infinite value.

\vskip5mm
{\bf Definition 9.1.} Let $(\Omega,\mathfrak M,\lambda)$ be a probability space. Given a measurable 
function $u:\Omega \rightarrow \R^n$, define its Luxemburg pseudo-norm by
\be
\|u\|_L = \|u\|_{L(\lambda)} = \inf\Big\{r > 0: \int L(u/r)\,d\lambda \leq 1\Big\},
\en
when the value of $r$ under the infimum sign exists, and otherwise $\|u\|_L = \infty$.
If $\|u\|_L$ is finite, define the Orlicz pseudo-norm by
\be
|u|_L = |u|_{L(\lambda)} = \sup\Big\{\int \left<u,v\right>d\lambda: \|v\|_{L^*} \leq 1\Big\}.
\en

\vskip2mm
The equalities (9.1)--(9.2) represent a natural generalization of the Luxemburg and Orlicz norms
of scalar functions and with $L$ being even. The scalar case is discussed in detail in 
Krasnoselskii and Rutickii \cite{K-R} and Maligranda \cite{M}. The latter book also contains
historical remarks: The definition (9.2) was first given in 1930's by Orlicz \cite{O1} assuming the
$\Delta_2$-condition and then in \cite{O2} for general Young functions.
The equivalent definition (9.1) was surprisingly given much later in 1950's in the work by
Nakano \cite{N} and Luxemburg \cite{Lu}.

As far we know, the multidimensional case was not considered in the literature so far.
This case, especially for non-smooth convex functions $L$, leads to some non-trivial questions
about subgradients of $L$ and the equivalence of the two pseudo-norms. Therefore,
in this section we collect several observations made in the recent paper \cite{B-G}, which was
a preliminary step towards the current work.

First let us note that in (9.2), $\|v\|_{L^*}$ is defined according to (9.1), although the condition 
$L^*(x)>0$ for all $x \neq 0$ may be violated.
If $\|u\|_L$ is finite and positive, the infimum in (9.1) is attained at $r = \|u\|_L$, so that
$$
\int L(u/\|u\|_L)\,d\lambda = 1 \quad (0 < \|u\|_L < \infty).
$$
It should be clear
that $\|u\|_L = 0$ or $|u|_L = 0$ implies that that $u=0$ with $\lambda$-probability 1. 

The integrals in (9.2) are well-defined and bounded from above when $\|u\|_L$ is finite. 
This restriction may be removed, since in this definition, 
one may additionally require that under the supremum sign the functions $v$ satisfy 
$\left<u,v\right> \geq 0$. That is,
\bee
|u|_L
 & = &
\sup\Big\{\int \left<u,v\right>d\lambda: \|v\|_{L^*} \leq 1,\, \left<u,v\right> \geq 0\Big\} \nonumber \\
 & = &
\sup\Big\{\int \left<u,v\right>^+ d\lambda: \|v\|_{L^*} \leq 1\Big\},
\ene
which may be taken as an equivalent definition instead of (9.2). Note that, due to the relation
$\left<x,y\right> \leq L(x) + L^*(y)$, $x,y \in \R^n$, we always have
$$
\int \left<u,v\right>^+  d\lambda \leq 2\,\|u\|_L \, \|v\|_{L^*}.
$$

Let us now list basic properties of the functionals in (9.1)-(9.2).

\vskip5mm
{\bf Proposition 9.2.} {\sl 
The functional $\|u\|_L$ is non-negative, positive homogeneous of order $1$ and convex:

\vskip2mm
$a)$ \ $0 \leq \|u\|_L \leq \infty$;

\vskip1mm
$b)$ \ $\|\alpha u\|_L = \alpha\, \|u\|_L$ for all $\alpha>0$;

\vskip1mm
$c)$ \ $\|t_1 u_1 + t_2 u_2\|_L  \leq t_1 \|u_1\|_L + t_2 \|u_2\|_L$ for all $t_1,t_2 \geq 0$ 
such that $t_1+t_2=1$.

\vskip1mm
$d)$\, Hence, this functional is subadditive: $\|u_1 + u_2\|_L \leq \|u_1\|_L + \|u_2\|_L$.

\vskip2mm
\noindent
The same properties hold for the functional $|u|_L$.
}

\vskip5mm
In this sense, $\|u\|_L$ and $|u|_L$ represent pseudo-norms. Since $L$ is not required to 
be even, i.e. $L(-x) = L(x)$, the property $\|-u\|_L = \|u\|_L$ may not hold. 
However, if $L$ is even, we obtain a norm in the Orlicz Banach space of all measurable 
functions $u$ on $\Omega$ such that $\|u\|_L < \infty$.

Restricting ourselves to dimension $n=1$, the most frequent choice of the convex  function $L$ 
is the power function $L(x) = |x|^p$, $1 < p < \infty$. Then the dual transform is given by
$$
L^*(x) = \frac{1}{q p^{q-1}}\, |x|^q, \quad q = \frac{p}{p-1}.
$$
Hence
$$
\|u\|_L = \|u\|_p = \Big(\int |u|^p\,d\lambda\Big)^{1/p}, \quad 
\|u\|_{L^*} = \frac{1}{p^{1/p}\, q^{1/q}}\,\|u\|_q,
$$
and
\be
|u|_L = p^{1/p}\, q^{1/q}\,\|u\|_p = p^{1/p}\, q^{1/q}\,\|u\|_L.
\en

The basic relationship between the two pseudo-norms is described in the following:

\vskip5mm
{\bf Proposition 9.3.} {\sl Suppose additionally that $L$ is everywhere finite.
For any measurable function $u:\Omega \rightarrow \R^n$, 
\be
\|u\|_L \leq |u|_L \leq 2\,\|u\|_L.
\en
}

\vskip0mm
By the one-dimensional example as in (9.3), the factor 2 cannot be improved
in (9.4). Indeed, the constant $p^{1/p}\, q^{1/q}$ is greater than 1
and takes the maximal value 2 for $p=q=2$.

A useful property of the Luxembourg norm is its quasi-concavity with respect to 
measures $\lambda$. 
Given a measurable function $u:\Omega \rightarrow \R^n$, consider the functional 
\be
S(\lambda) = \|u\|_{L(\lambda)}
\en 
on the space of all probability measures $\lambda$ on $(\Omega,\frak M)$.

\vskip5mm
{\bf Proposition 9.4.} {\sl The functional $S$ is quasi-concave: Given probability measures 
$\lambda_1,\lambda_2$ on $\Omega$, for all $t_1,t_2 \geq 0$ such that $t_1 + t_2  = 1$, 
\be
S(t_1 \lambda_1 + t_2 \lambda_2) \geq \min\{S(\lambda_1),S(\lambda_2)\}.
\en
Moreover, if $L$ is positive homogeneous of order $p \geq 1$,
this functional is concave:
\be
S(t_1 \lambda_1 + t_2 \lambda_2) \geq t_1 S(\lambda_1) + t_2 S(\lambda_2).
\en
}

\vskip0mm
The inequality (9.6) can be extended, as well as (9.7), to arbitrary (``continuous") convex mixtures 
of probability measures. Let us state such a relation for convolutions $\lambda * \kappa$ on 
the Euclidean space $\Omega = \R^m$, 
which we equip with the Borel $\sigma$-algebra. Note that any such convolution represents 
a convex mixture of shifts or translates of $\lambda$ with a mixing measure $\kappa$.

\vskip5mm
{\bf Proposition 9.5.} {\sl Suppose that $L$ is positive homogeneous of order $p \geq 1$.
Given a Borel measurable function $u:\R^m \rightarrow \R^n$,
for all probability measures $\lambda$ and $\kappa$ on $\R^m$,
\be
\|u\|_{L(\lambda * \kappa)} \, \geq \, \int \|u(x-y)\|_{L(\lambda(dx))}\,d\kappa(y).
\en
}

\vskip0mm
Under the two-sided $\Delta_2$-condition, the relations (9.6)-(9.7) are extended in 
a somewhat weaker form. Given $p_1 \geq p_0 \geq 1$, define the quantity
\be
\gamma(p_1,p_0) = \sup \big\{a + b: a+br \leq \min(r^{p_1},r^{p_0}) \ {\rm for \ all} \ r \geq 0\big\}.
\en
In particular, $0 < \gamma(p_1,p_0) \leq 1$ and $\gamma(p_1,p_0) = 1$ if $p_1 = p_0$ (and only in this case).
We use this definition with
$$
p_1 = p_+ = \Phi_L'(1+), \quad
p_0 = p_- = \Phi_L'(1-).
$$

\vskip5mm
{\bf Proposition 9.6.} {\sl Given probability measures $\lambda_1,\dots,\lambda_k$ on 
$\Omega$, for all $t_i \geq 0$ such that $t_1 + \dots + t_k  = 1$, we have
$$
S(t_1 \lambda_1 + \dots + t_k \lambda_k) \geq \gamma \sum_{i=1}^k t_i S(\lambda_i)
$$
with constant $\gamma = \gamma(p_+,p_-)$. As a consequence, if $u:\R^m \rightarrow \R^n$ is
Borel measurable, then for all probability measures $\lambda$ and $\kappa$ on $\R^m$,
\be
\|u\|_{L(\lambda * \kappa)} \, \geq \, \gamma \int \|u(x-y)\|_{L(\lambda(dx))}\,d\kappa(y).
\en
}

\vskip7mm
\section{{\bf The Dual Sobolev Norms $\|\nu - \mu\|_{H^{-1,p}(\lambda)}$}}
\setcounter{equation}{0}

\vskip2mm
\noindent
Let $L:\R^n \rightarrow [0,\infty)$ be a convex function such that $L(0)=0$,
$L(x)>0$ for all $x \neq 0$, and let
$$
L^*(x) = \sup_{y \in \R^n} \big[\left<x,y\right> - L(y)\big], \quad x \in \R^n.
$$

\vskip2mm
{\bf Definition 10.1.} Given three probability measures $\mu,\nu,\lambda$ on $\R^n$, define 
the dual Sobolev  pseudo-norm
\be
\|\nu - \mu\|_{H^{-1,L}(\lambda)} = \sup\Big\{\int f\,d\nu - \int f\,d\mu:
\int L^*(\nabla f)\,d\lambda \leq 1\Big\},
\en
where the supremum is running over all 
$C^\infty$-smooth functions $f$ on $\R^n$ with compact support subject to the above integral condition.

\vskip5mm
This definition is a bit more general in comparison with (1.3) where $\lambda = \mu$.
We consider this more general setting with regard to some other applications.

Thus, using the Luxemburg pseudo-norm generated by $L^*$,
the value $c = \|\nu - \mu\|_{H^{-1,L}(\lambda)}$ represents an optimal constant in the inequality
\be
\int f\,d\nu - \int f\,d\mu \leq c\,\|\nabla f\|_{L^*(\lambda)}
\en 
in the class of all $C^\infty$-smooth functions $f$ on $\R^n$ with compact support.
In this section, we assume that $L^*$ is everywhere finite.

\vskip5mm
{\bf Proposition 10.2.} {\sl The inequality $(10.2)$ may be extended with the same constant $c$ to
the class of all bounded Lipschitz $C^1$-smooth functions $f$ on $\R^n$. Moreover,
if the measure $\lambda$ is absolutely continuous, $(10.2)$ holds true for
all bounded Lipschitz functions $f$ on $\R^n$.
}

\vskip5mm
{\bf Proof.} 
First, we derive (10.2) for any $C^1$-smooth function $f$ on $\R^n$ 
with compact support. In view of the homogeneity of this inequality, we may assume that
$\|\nabla f\|_{L^*(\lambda)} = 1$.

Take a $C^\infty$-smooth function $\omega:\R^n \rightarrow [0,\infty)$ supported on the unit ball 
$B_1$ and such that $\int \omega(x)\,dx = 1$. The functions 
$\omega_\ep(x) = \ep^n \omega(\ep x)$, $\ep > 0$, are supported 
on the balls  $B_\ep = \{x \in \R^n: |x| \leq \ep\}$. Hence, the  convolutions
\be
f_\ep(x) = (f * \omega_\ep)(x) = \int f(x-z)\,\omega_\ep(z)\,dz = \int \omega_\ep(x-z)\,f(z)\,dz
\en
are also $C^\infty$-smooth and have a compact support. Thus, 
according to the definition,
\be
\int f_\ep\,d\nu - \int f_\ep\,d\mu \leq c\,\|\nabla f_\ep\|_{L^*(\lambda)}.
\en 

By the continuity of $f$, we have $f_\ep(x) \rightarrow f(x)$ as $\ep \rightarrow 0$ 
and $\sup_\ep |f_\ep(x)| \leq M$ for all $x \in \R^n$ with some constant $M>0$. 
Hence, by the dominated convergence theorem,
$$
\lim_{\ep \rightarrow 0}\int f_\ep\, d\mu = \int f d\mu, \quad
\lim_{\ep \rightarrow 0}\int f_\ep\, d\nu = \int f d\nu.
$$

Using the Lipschitz property, one may also differentiate in (10.3) to get a similar representation
$$
\nabla f_\ep(x) = \int \nabla f(x-z)\,\omega_\ep(z)\,dz.
$$
By the continuity of the gradient $\nabla f$, we have $\nabla f_\ep(x) \rightarrow \nabla f(x)$
and $L^*(\nabla f_\ep(x)) \rightarrow L^*(\nabla f(x))$ as $\ep \rightarrow 0$ 
for any $x \in \R^n$. Since the gradient is bounded, the dominated convergence theorem is applicable, 
and we get
$$
\lim_{\ep \rightarrow 0}\int L^*(\nabla f_\ep)\, d\lambda = \int L^*(\nabla f)\, d\lambda = 1.
$$
Given $\delta>0$, we conclude that the integrals on the above left-hand side are smaller 
than $1+\delta$ for all $\ep>0$ small enough. But then, by the convexity of $L^*$,
$\int L^*(\nabla f_\ep/(1+\delta))\, d\lambda \leq 1$, which means that
$\|\nabla f_\ep\|_{L^*(\lambda)} \leq 1+\delta$. In other words,
$$
\limsup_{\ep \rightarrow 0} \, \|\nabla f_\ep\|_{L^*(\lambda)} \leq 1+\delta
$$
and hence this lim\,sup is at most 1. Returning to (10.4), we thus obtain (10.2) for $f$.

Now, given a bounded Lipschitz $C^1$-smooth function $f$ on $\R^n$, consider the functions
$f_r(x) = f(x)\,g(x/r)$, $x \in \R^n$, $r>0$, where $g:\R^n \rightarrow [0,1]$ is $C^1$-smooth,
supported on the unit ball, and is such that $g(x) = 1$ for $|x| \leq 1/2$. By the previous step,
we have (10.4) for all $f_r$. They are supported on the balls $B_r$ and Lipschitz.
In addition,  $|\nabla f_r(x)| \leq |\nabla f(x)| + M/r$ for all $x \in \R^n$ and $r>0$ with some
constant $M$. Then letting $r \rightarrow \infty$, we obtain the same relation for $f$.

Turning to the last claim, suppose that $\lambda$ is absolutely continuous with respect to the
Lebesgue measure on $\R^n$. By the Rademacher theorem, Lipschitz functions are differentiable 
almost everywhere, so, the Luxemburg norm in (10.2) is well-defined. One can repeat the previous 
argument with the following refinement: If the smoothing density $\omega = \omega(|x|)$ 
depends on $|x|$, then $\nabla f_\ep(x) \rightarrow \nabla f(x)$ and therefore 
$L^*(\nabla f_\ep(x)) \rightarrow L^*(\nabla f(x))$ as $\ep \rightarrow 0$ for almost all $x \in \R^n$. 
\qed

\vskip2mm
The following property describes the lower semi-continuity of the dual Sobolev norm
in the topology of weak convergence of probability measures on $\R^n$.

\vskip5mm
{\bf Proposition 10.3.} {\sl Suppose that three sequences $(\mu_k)_{k \geq 1}$, 
$(\nu_k)_{k \geq 1}$,  $(\lambda_k)_{k \geq 1}$ of probability measures on $\R^n$ are weakly 
convergent as $k \rightarrow \infty$ to probability measures $\mu$, $\nu$, $\lambda$. Then, 
\be
\|\nu - \mu\|_{H^{-1,L}(\lambda)} \leq \liminf_{k \rightarrow \infty}\, 
\|\nu_k - \mu_k\|_{H^{-1,L}(\lambda_k)}.
\en
}

\vskip0mm
{\bf Proof.} Let $c = \|\nu - \mu\|_{H^{-1,L}(\lambda)}$ be positive and $0 < c' < c$. 
By Definition 10.1, we may choose a $C^\infty$-smooth function $f$ on $\R^n$ with compact 
support such that $\|\nabla f\|_{L^*(\lambda)} = 1$ and $\int f\,d(\nu - \mu) > c'$. 
Since the functions $f$ and $L^*(\nabla f)$ 
are bounded and continuous, the weak convergence implies that as $k \rightarrow \infty$
$$
\int f\,d(\nu_k - \mu_k) \rightarrow \int f\,d(\nu - \mu) \quad {\rm and} \quad
\int L^*(\nabla f)\,d\lambda_k \rightarrow \int L^*(\nabla f)\,d\lambda \, = \, 1.
$$
Hence, for any $\ep > 0$, there is an integer $N \geq 1$ such that
$\int L^*(\nabla f)\,d\lambda_k \leq 1+\ep$  for all $k \geq N$ and also $\int f\,d(\nu_k - \mu_k) \geq c'$.
Using $L^*(\alpha x) \leq \alpha L^*(x)$ which holds by the convexity of $L^*$ for all 
$x \in \R^n$ and $\alpha \in [0,1]$, we get
$$
\int L^*\Big(\frac{\nabla f}{1+\ep}\Big)\,d\lambda_k \leq 1.
$$
This means that $\|\nabla f\|_{L^*(\lambda_k)} \leq 1+\ep$, and as a consequence, for all $k \geq N$,
$$
\|\nu_k - \mu_k\|_{H^{-1,L}(\lambda_k)} \geq 
\frac{1}{\|\nabla f\|_{L^*(\lambda_k)}}\,\int f\,d(\nu_k - \mu_k) > \frac{c'}{1+\ep}.
$$
Hence
$$
\liminf_{k \rightarrow \infty}\, \|\nu_k - \mu_k\|_{H^{-1,L}(\lambda_k)} \geq \frac{c'}{1+\ep}.
$$
Letting $\ep \rightarrow 0$ and $c' \rightarrow c$, we arrive at (10.5).
\qed

\vskip7mm
\section{{\bf Continuity of the Dual Sobolev Norm along Convolutions}}
\setcounter{equation}{0}

\vskip2mm
\noindent
For convolutions of  probability measures  $\mu$, $\nu$, $\lambda$ on $\R^n$, the inequality 
(10.5) may be reversed. More precisely, we fix a probability measure $\kappa$ on $\R^n$ 
and denote by $\kappa_\ep$ the image of $\kappa$ under the linear map
$x \rightarrow \ep x$ with parameter $\ep>0$. Define the convolutions
$$
\mu_\ep = \mu * \kappa_\ep, \quad \nu_\ep = \nu * \kappa_\ep, \quad 
\lambda_\ep = \lambda * \kappa_\ep.
$$
First let us start with the homogeneous case.

\vskip5mm
{\bf Proposition 11.1.} {\sl If $L$ is a positive homogeneous convex function on $\R^n$ 
of order $p > 1$, then
\be
\|\nu - \mu\|_{H^{-1,L}(\lambda)} = \lim_{\ep \rightarrow 0}\, 
\|\nu_\ep - \mu_\ep\|_{H^{-1,L}(\lambda_\ep)}.
\en
}

\vskip0mm
{\bf Proof.} Since the convolved measures are weakly convergent as $\ep \rightarrow 0$
to $\mu$, $\nu$, $\lambda$, we have, by Proposition 10.3,
\be
\|\nu - \mu\|_{H^{-1,L}(\lambda)} \, \leq \, \liminf_{\ep \rightarrow 0}\, 
\|\nu_\ep - \mu_\ep\|_{H^{-1,L}(\lambda_\ep)}.
\en

For the opposite direction, take a $C^\infty$-smooth function $f$ on $\R^n$ with compact support
and apply (10.2) to the function $x \rightarrow f(x-y)$. It gives
$$
\int f(x-y)\,d\nu(x) - \int f(x-y)\,d\mu(x) \, \leq \, c\,\|\nabla f(x-y)\|_{L^*(\lambda(dx))}
$$
with constant $c = \|\nu - \mu\|_{H^{-1,L}(\lambda)}$.
Integration of this inequality over the measure $\kappa_\ep$ yields
\be
\int f\,d\nu_\ep - \int  f\,d\mu_\ep \, \leq \, 
c \int \|\nabla f(x-y)\|_{L^*(\lambda(dx))}\,d\kappa_\ep(y).
\en

Recall that $L^*$ is positive homogeneous convex function on $\R^n$ of order $q = p/(p-1)$,
according to Proposition 7.5. Now we may use the concavity of the functional $S$ in (9.5) and
apply Proposition 9.5 with $u = \nabla f$
and with $\kappa_\ep$ in place of $\kappa$. By (9.8), for any $\ep>0$,
$$
\int \|\nabla f(x-y)\|_{L^*(\lambda(dx))}\,d\kappa_\ep(y) \leq 
\|\nabla f\|_{L^*(\lambda * \kappa_\ep)}.
$$
Inserting this bound in (10.3), we get
$$
\int f\,d\nu_\ep - \int  f\,d\mu_\ep \, \leq \, c\,\|\nabla f\|_{L^*(\lambda _\ep)}
$$
Since $f$ was arbitrary, this bound amounts to the upper bound
$\|\nu_\ep - \mu_\ep\|_{H^{-1,L}(\lambda_\ep)} \leq c$. Hence 
$$
\limsup_{\ep \rightarrow 0}\, \|\nu_\ep - \mu_\ep\|_{H^{-1,L}(\lambda_\ep)} \leq c.
$$
Together with (11.2) this proves (11.1).
\qed

\vskip5mm
In a weaker form, the property (11.1) can be extended to a larger class of convex functions $L$ 
satisfying the two-sided $\Delta_2$-condition. Recall that in this case we considered the Young 
function $\Phi_L(r)$ and its one-sided derivatives $p_+ = \Phi_L'(1+)$ and $p_- = \Phi_L'(1-)$.

\vskip5mm
{\bf Proposition 11.2.} {\sl Let $L$ be a non-negative convex function on $\R^n$ such that
$L(0) = 0$ and $L(x) > 0$, satisfying the two-sided $\Delta_2$-condition with $p_-> 1$. Then
\begin{eqnarray}
\|\nu - \mu\|_{H^{-1,L}(\lambda)} 
 & \leq & 
\liminf_{\ep \rightarrow 0}\ \|\nu_\ep - \mu_\ep\|_{H^{-1,L}(\lambda_\ep)} \nonumber \\
& \leq &
\limsup_{\ep \rightarrow 0}\, \|\nu_\ep - \mu_\ep\|_{H^{-1,L}(\lambda_\ep)} \, \leq \, 
\frac{1}{\gamma}\, \|\nu - \mu\|_{H^{-1,L}(\lambda)},
\end{eqnarray}
where the constant $\gamma = \gamma(p_+,p_-)$ was defined in $(9.9)$.
}

\vskip5mm
{\bf Proof.} By Proposition 7.2, $L(rx) \geq r^{p_-} L(x)$ for all $r \geq 1$ and 
$x \in \R^n$, implying that $L$ satisfies the super-linear growth condition. Then $L^*$ is finite
everywhere, so that we may apply Proposition 10.3. Since the convolved measures are weakly 
convergent as $\ep \rightarrow 0$ to $\mu$, $\nu$ and $\lambda$, we thus obtain
the first inequality in (11.4).

To derive the last inequality, one may repeat the previous arguments used in the proof
of Proposition 11.1. First we have the upper bound (11.3). Secondly, by the inequality (9.10)
of Proposition 9.6, for any $\ep>0$,
$$
\int\!\!\!\int \|\nabla f(x-y)\|_{L^*(\lambda(dx))}\,d\kappa_\ep(y) \leq \frac{1}{\gamma}\,
\|\nabla f\|_{L^*(\lambda * \kappa_\ep)}.
$$
Inserting this bound in (11.3), we get
$$
\int f\,d\nu_\ep - \int  f\,d\mu_\ep \, \leq \, \frac{c}{\gamma}\,\|\nabla f\|_{L^*(\lambda _\ep)}.
$$
Since $f$ was arbitrary, this bound amounts to the upper bound
$\|\nu_\ep - \mu_\ep\|_{H^{-1,L}(\lambda_\ep)} \leq \frac{c}{\gamma}$. Hence 
$$
\limsup_{\ep \rightarrow 0}\, \|\nu_\ep - \mu_\ep\|_{H^{-1,L}(\lambda_\ep)} \leq \frac{c}{\gamma},
$$
thus proving the last inequality in (11.4).
\qed

\vskip7mm
\section{{\bf Proof of Theorem 1.2 and Theorem 1.3}}
\setcounter{equation}{0}

\vskip2mm
{\bf Proof of Theorem 1.2.}
To derive (1.4), we may assume that $c = \|\nu - \mu\|_{H^{-1,L}(\mu)}$ is finite.
Recall the basic dual representation
\be
\T_L(\mu,\nu) = \sup \bigg[\int Q_1 f\,d\nu - \int f\,d\mu\bigg]
\en
as in Proposition 3.1, where $Q_1$ is the infimum-convolution operator associated 
to $L$ at time $t=1$, and where the supremum is running over all bounded continuous functions 
$f$ on $\R^n$ (equivalently, over all $C^\infty$-smooth non-negative functions with compact support).

To bound the difference of integrals in (12.1), we adapt M. Ledoux' argument, although
this has required a significant development of the theory of Luxemburg and Orlicz pseudo-norms
for vector-valued functions.

Introduce the function $g = \varphi - 1$, where $\varphi = d\nu/d\mu$ is a density of $\nu$ 
with respect to $\mu$ and note that $g \geq -1$ $\mu$-a.e. with $\int g\,d\mu = 0$. 

Let $f$ be a bounded Lipschitz function on $\R^n$ and let $\theta:[0,1] \rightarrow [0,1]$ 
be a continuous increasing function, with a bounded continuous derivative $\theta'(t)$ in 
$0 < t < 1$, and such that $\theta(0) = 0$, $\theta(1)=1$.
For a fixed $x \in \R^n$, consider the function
$$
u_x(t) = \big(1 + \theta(t)\, g(x)\big)\, Q_t f(x) - f(x).
$$
It is Lipschitz on $[0,1]$. By the Hamilton-Jacobi equation (Proposition 6.1), it has derivative 
$$
u_x'(t) \, = \,
\theta'(t)\, g(x)\, Q_t f(x) - \big(1 + \theta(t)\, g(x)\big) L^*(\nabla Q_t f(x))
$$
for all $t>0$ such that the function $(x,t) \rightarrow Q_t f(x)$ is differentiable at the point $(x,t)$.
By Corollary 6.2, for almost all $x \in \R^n$, 
$u_x'(t)$ also serves as a Radon-Nikodym derivative of $u_x$. Since $u_x(0) = 0$ and 
$$
u_x(1) = \varphi(x) Q_1 f(x) - f(x), 
$$
we may write $u_x(1) = \int_0^1 u_x'(t)\,dt$, that is,
$$
\varphi(x) Q_1 f(x) - f(x) \, = \, g(x) \int_0^1 \theta'(t) Q_t f(x)\,dt -
\int_0^1 (1 + \theta(t) g)\,L^*(\nabla Q_t f(x))\,dt,
$$
which holds true for almost all $x \in \R^n$. Here $Q_t f(x)$ and $|\nabla Q_t f(x)|$ 
are uniformly bounded over all $(x,t)$ and Borel measurable. Integrating this equality 
over $\mu$ and using the absolute continuity of this measure, we get
\begin{eqnarray}
I(f) \ \equiv \ \int Q_1 f\,d\nu - \int f\,d\mu
 & = &
\int\!\!\! \int_0^1 \theta'(t)\, g(x)\, Q_t f(x)\,dt\,d\mu(x) \nonumber \\
 & & \hskip-15mm - \
\int\!\!\! \int_0^1 \big(1 + \theta(t)\, g(x)\big) L^*(\nabla Q_t f(x))\,dt\,d\mu(x).
\end{eqnarray}

Note that the second last integral in (12.2) with respect to $\mu$ is equal to and satisfies
$$
\int g(x)\, Q_t f(x)\,d\mu(x) \, = \, \int Q_t f\,d\nu - \int Q_t f\,d\mu  \, \leq \,
c\,\|\nabla Q_t f\|_{L^*(\mu)},
$$
where $\|\nabla Q_t f\|_{L^*(\mu)}$ is the Luxemburg pseudo-norm over the measure $\mu$
corresponding to the convex function $L^*$. In this step we applied the last claim
in Proposition 10.2, where the assumption that the measure $\lambda = \mu$ is absolutely
continuous was needed. Hence, using also  $1 + \theta(t)\, g(x) \geq 1 - \theta(t) > 0$ a.e. in the last 
integrand in (12.2), this inequality yields
\begin{eqnarray}
I(f) 
 & \leq & 
c \int_0^1 \theta'(t)\, \|\nabla Q_t f\|_{L^*(\mu)}\,dt \nonumber \\
 & & - \
\int_0^1 \big(1 - \theta(t)\big)\,dt \int  L^*(\nabla Q_t f(x))\,d\mu(x).
\end{eqnarray}

Let us rewrite this inequality as $I \, \leq \, \int_0^1 I_t(f)\,dt$ in terms of the functionals
$$
I_t(f) \, = \, c\theta'(t)\, \|\nabla Q_t f\|_{L^*(\mu)} - 
\big(1 - \theta(t)\big) \int L^*(\nabla Q_t f)\,d\mu.
$$
Using $L^{**} = L$, we relate the Luxemburg pseudo-norm to the Orlicz pseudo-norm
$$
|\nabla Q_t f|_{L^*(\mu)} = \sup \Big\{\int \left<\nabla Q_t f,v\right> d\mu:
\|v\|_{L(\mu)} \leq 1 \Big\},
$$
where the supremum is running over all Borel measurable maps $v:\R^n \rightarrow \R^n$
with $\|v\|_{L(\mu)} \leq 1$. Namely, according to Proposition 9.3,
$$
\|\nabla Q_t f\|_{L^*(\mu)} \leq |\nabla Q_t f|_{L^*(\mu)}. 
$$
Therefore, 
\begin{eqnarray}
I_t(f) 
 & \leq &
c \theta'(t)\sup_{\|v\|_{L(\mu)} \leq 1} 
\int \left<\nabla Q_t f,v\right> d\mu - (1 - \theta(t)) \int L^*(\nabla Q_t f)\,d\mu \nonumber \\
 & & \hskip-20mm = \
(1-\theta(t)) \sup_{\|v\|_{L(\mu)} \leq 1} \,
\int \Big(\frac{c \theta'(t)}{1 - \theta(t)} \left<\nabla Q_t f(x),v(x)\right> - 
L^*(\nabla Q_t f(x))\Big)\, d\mu(x).
\end{eqnarray}
By the definition of the Legendre transform and using $L^{**} = L$ once more, we have
$$
a \left<u,v\right> - L^*(u) \leq L(av), \quad a \in \R, \ u,v \in \R^n.
$$
Applying this with $u = \nabla Q_t f(x)$, $v = v(x)$ for an arbitrary fixed $x \in \R^n$,
and with constant $a = c \theta'(t)/(1 - \theta(t))$, we can bound the integrand in (12.4), by
\be
L\Big(\frac{c \theta'(t)}{1 - \theta(t)}\,v(x)\Big) \, \leq \, 
\Phi\Big(\frac{c \theta'(t)}{1 - \theta(t)}\Big) L(v(x)),
\en
where $\Phi = \Phi_L$. Note that, under the supremum sign in (12.4), the condition on $v$
means that $\int L(v)\,d\mu(x) \leq 1$. Thus, after integration over $\mu$ and using (12.5),
(12.4) yields
\be
I_t(f) \leq (1-\theta(t))\, \Phi\Big(\frac{c \theta'(t)}{1 - \theta(t)}\Big), \quad 0 < t < 1.
\en

Recalling (12.3), further integration of (12.6) over $t$ in the unit interval gives
$$
\int Q_1 f\,d\nu - \int f\,d\mu \leq 
\int_0^1 (1 - \theta(t)) \, \Phi\Big(\frac{c \theta'(t)}{1 - \theta(t)}\Big)\,dt.
$$
Taking the supremum over all admissible $f$ and recalling (12.1), we arrive at
$$
\T_L(\mu,\nu) \leq 
\int_0^1 (1 - \theta(t)) \, \Phi\Big(\frac{c \theta'(t)}{1 - \theta(t)}\Big)\,dt.
$$
It remains to apply Propositions 8.2--8.3: For the special function $\theta(t)$ the above integrand
is a constant which does not exceed $A_L \Phi(c)$.
\qed

\vskip5mm
{\bf Proof of Theorem 1.2.}
Let $\kappa_\ep$ denote the Gaussian measure on $\R^n$ with mean zero and
covariance matrix $\ep^2\,{\rm Id}$, $\ep>0$. Consider the convolutions
$\mu_\ep = \mu * \kappa_\ep$ and $\nu_\ep = \nu * \ep$. Since $L$ has at most a polynomial growth,
the hypothesis $\int\!\!\int L(x-y)\,d\mu_\ep(x)\,d\nu_\ep(y) < \infty$ is satisfied. 

We are in position to apply Theorem 1.2 to the couple $(\mu_\ep,\nu_\ep)$, which gives
\be
\T_L(\mu_\ep,\nu_\ep) \leq \Phi(p_+)\, \Phi_L\big(\|\nu_\ep - \mu_\ep\|_{H^{-1,L}(\mu_\ep)}\big).
\en
It should be clear that
\be
\T_L(\mu_\ep,\nu_\ep) \rightarrow \T_L(\mu,\nu) \quad {\rm as} \ \ep \rightarrow 0.
\en
On the other hand, by the upper bound in (11.4),
$$
\limsup_{\ep \rightarrow 0}\, \|\nu_\ep - \mu_\ep\|_{H^{-1,L}(\mu_\ep)} \, \leq \, 
\frac{1}{\gamma}\, \|\nu - \mu\|_{H^{-1,L}(\mu)},
$$
implying
\begin{eqnarray}
\limsup_{\ep \rightarrow 0}\, \Phi_L(\|\nu_\ep - \mu_\ep\|_{H^{-1,L}(\mu_\ep)}) 
 & \leq &
\Phi_L\Big(\frac{1}{\gamma}\, \|\nu - \mu\|_{H^{-1,L}(\mu)}\Big) \nonumber \\
 & \leq &
\Phi_L(1/\gamma)\, \|\nu - \mu\|_{H^{-1,L}(\mu)}.
\end{eqnarray}
In the last step, we used the sub-multiplicativity of the function $\Phi_L$.
It remains to apply (12.8)--(12.9) in (12.7) and to turn to the limit as $\ep \rightarrow 0$.
\qed

\end{document}